\documentstyle[amscd,amssymb,verbatim,11pt]{amsart}

\theoremstyle{plain}
\newtheorem{Prop}{Proposition}[section]
\newtheorem{Thm}[Prop]{Theorem}
\newtheorem{Cor}[Prop]{Corollary}

\newtheorem{Lem}[Prop]{Lemma}

\theoremstyle{definition}
\newtheorem{Def}[Prop]{Definition}

\theoremstyle{remark}
\newtheorem{Rem}[Prop]{Remark}
\newtheorem{ex}[Prop]{Example}

\def\asdim{\text{asdim}}
\def\cdim{\text{cdim}}
\def\dim{\text{dim}}
\def\diam{\text{diam}}

\def\NN{{\mathbb N}}
\def\QQ{{\mathbb Q}}
\def\RR{{\mathbb R}}
\def\la{{\lambda}}
\def\La{{\Lambda}}

\def\UU{{\mathcal U}}

\def\ZZ{{\mathbb Z}}

\begin{document}
\title[Assouad-Nagata dimension of connected Lie groups]{Assouad-Nagata dimension of connected Lie groups}

\author{J.Higes}
\email{josemhiges@@yahoo.es}

\author{I. Peng}
\email{ipeng@@indiana.edu}



\thanks{The first named author is supported by project MEC, MTM2006-0825 and 'contrato flechado' i-math.}

\begin{abstract} We prove that the asymptotic Assouad-Nagata dimension of a connected  Lie group $G$ equipped with a left-invariant Riemannian metric coincides
with its topological dimension of $G/C$ where $C$ is a maximal compact subgroup. To prove it we will compute the Assouad-Nagata dimension of connected solvable Lie groups and semisimple Lie groups. As a consequence we show that the asymptotic Assouad-Nagata dimension of a polycyclic group equipped with a word
metric is equal to its Hirsch length and that some wreath-type finitely generated groups can not be quasi-isometric to any cocompact lattice on a connected Lie group.
\end{abstract}

\maketitle

\tableofcontents
\section{Introduction}
The Assouad-Nagata dimension was introduced by Assouad in \cite{Assouad} inspired  from
the ideas of Nagata. Metric spaces of finite Assouad-Nagata dimension satisfy interesting
geometric properties. For example they admit quasisymmetric embeddings into the product
of finitely many trees \cite{Lang} and have nice Lipschitz extension properties (see
\cite{Lang} and \cite{Brod-Dydak-Higes-Mitra}). The class of metric spaces with finite
Nagata dimension includes in particular all doubling spaces, metric trees, Euclidean
buildings, and homogeneous or pinched negatively curved Hadamard manifolds. \par In last
years a small scale version and a large scale version of the Assouad-Nagata dimension
have been the focus of interesting research. The small scale version has been studied in
the framework of hyperbolic groups under the name of capacity dimension (see
\cite{Buyalo2}). The large scale version of Assouad-Nagata dimension has been referred to
by several names, such as asymptotic dimension of linear type, asymptotic dimension with
Higson property or asymptotic Assouad-Nagata dimension. This last name would be the one
we will use in this paper. \par The asymptotic Assouad-Nagata dimension is a linear
version of the asymptotic dimension, and it is invariant under quasi-isometries. But
while the asymptotic dimension remains invariant under coarse equivalences, the
asymptotic Assouad-Nagata dimension does not. Therefore it is reasonable to expect that
there would be more relationships between the asymptotic Assouad-Nagata dimension and
other quasi-isometric invariants of geometric group theory.  For example in \cite{Nowak}
and \cite{Brod-Dydak-Lang} the asymptotic Assouad-Nagata dimension was related with the
growth type of amenable groups and wreath products, in \cite{Dydak-Higes} it was shown
that the asymptotic Assouad-Nagata dimension bounds from above the topological dimension
of the asymptotic cones.  From the results of \cite{Gal} it is easy to see that every
metric space of finite asymptotic Assoaud-Nagata dimension has Hilbert compression equals
one. \par On the other hand, estimating the asymptotic Assouad-Nagata dimension of a
group is more difficult. Most of the methods developed to calculate the asymptotic
dimension of a group behaves poorly for asymptotic Assouad-Nagata dimension. For example,
in \cite{BD} Bell and Dranishnikov developed a technique to estimate the asymptotic
dimension of the fundamental group of a graph of groups. Unfortunately such method can
not be applied directly to the asymptotic Assouad-Nagata dimension. \par For connected
Lie  groups, it is natural to study short exact sequences of the form:
$$1 \to H \to G \to G/H \to 1.$$  The idea is to estimate the dimension of $G$ from the dimensions of $H$ and $G/H$. For example if $G$ is nilpotent
then $H$ could be abelian and $G/H$ nilpotent with lower degree of nilpotency. \par In
\cite{Brod-Dydak-Levin-Mitra}, Brodskyi, Dydak, Levin and Mitra developed techniques to
study the asymptotic Assouad-Nagata dimension of short exact sequences of finitely
generated groups with word metrics, and the ideas were applied successfully in
\cite{Dydak-Higes} to calculate the asymptotic Assouad-Nagata dimension of the Heisenberg
group.  But the main difficulty in the general case is to understand the distortion of
the subgroup $H$ in $G$.  \par In this paper we study the Assouad-Nagata dimension (at
small and large scale) of connected Lie groups equipped with left invariant Riemannian
metrics. First we will analyze the Assouad-Nagata dimension of simply connected solvable
Lie groups. The key tool for such goal is a generalization of the results from
\cite{Brod-Dydak-Levin-Mitra} in the setting of finitely generated groups with word
metrics, to general topological groups with left invariant metrics. Using this and some
facts from differential geometry, we will show that the Assouad-Nagata dimension of a
connected solvable Lie group is equal to its topological dimension.  As a consequence of
this we will prove that the asymptotic Assouad-Nagata dimension of a polycyclic group is
equal to its Hirsch length. This answers an open problem of the asymptotic Assouad-Nagata
dimension  as Question 4 of \cite{Nowak} and Problem 8.3 of \cite{DranProblems}(notice
that in \cite{DranProblems} it is said that Osin solved such problem but no proof has
been provided). Also our results increase the catalogue of finitely generated groups with
finite asymptotic Assouad-Nagata dimension. So far, the classes of groups known to have
finite asymptotic Assouad-Nagata dimension are Coexeter groups, abelian groups,
hyperbolic groups, free groups and some types of Baumslag-Solitar groups.  Recently it
was shown that the asymptotic Assouad-Nagata dimension is preserved under free products
(see \cite{Brod-Higes}). \par The next step in the proof of the main theorem of this
paper will be the study of the Assouad-Nagata dimension of semisimple Lie groups. For
such goal we will use the Iwasawa decomposition of a semisimple Lie group. Roughly
speaking an Iwasawa decomposition will say that a semisimple Lie group is equivalent to
the product of a connected solvable group and a special group $K$ for which we can apply
the methods of \cite{Brod-Dydak-Levin-Mitra}.\par

After studying the Assouad-Nagata dimension of semisimple Lie groups we will focus in the
main result of this paper: the Assouad-Nagata dimension of connected Lie groups. Notice
that in \cite{Car}, Carlsson and Goldfarb proved that the asymptotic dimension of a Lie
group $G$ is equal to the topological dimension of $G/C$ where $C$ is its maximal compact
subgroup. Therefore it is natural to see if the results of \cite{Car} can be extended to
the asymptotic Assouad-Nagata dimension.  Unfortunately the same techniques of \cite{Car}
cannot be applied directly for asymptotic Assouad-Nagata dimension. The main tools of
\cite{Car} use strongly the invariance under coarse equivalences of the asymptotic
dimension. For example Proposition 3.4. of \cite{Car} uses such  property. We will show
in Example \ref{ExampleCarlsson} that also the techniques from Theorem 3.5. of \cite{Car}
can not be applied.\par As it was mentioned before, the asymptotic Assouad-Nagata
dimension has relationships with the topological dimension of the asymptotic cones. In
\cite{Dydak-Higes} it was shown that the topological dimension of the asymptotic cone is
bounded from above by the asymptotic Assouad-Nagata dimension of the group. De Cornulier
computed in \cite{DeCor} the topological dimension of the asymptotic cones of connected
Lie groups. He defined the exponential radical of a connected Lie group and used it to
compute the dimension of the cones.  We will also use the exponential radical to compute
the Assouad-Nagata dimension of connected Lie groups. For this final theorem the
computations made for connected solvable Lie groups and semisimple Lie groups will play
an important role. \par  We will prove that the Assouad-Nagata dimension for a connected
Lie group $G$ is equal to the topological dimension of $G/C$, where $C$ is a maximal
compact subgroup. Hence our results can  be seen as a bridge among the ones of
\cite{Car}, \cite{DeCor} and \cite{Dydak-Higes}. Moreover it is not difficult to
show (although tedious) that combining the techniques of \cite{DeCor} and
\cite{Dydak-Higes} with our main theorem we can improve slightly the results of
\cite{DeCor} to compute the Assouad-Nagata dimension of the asymptotic cones of connected
Lie groups. \par

In Section \ref{Nagata}  we will study the behaviour of the Assouad-Nagata dimension
under some transformations on the metric.  Such results will be important when we study
the restriction of a  metric of a group to one of its subgroups. In Section
\ref{Hurewicz} we will generalize the results of \cite{Brod-Dydak-Levin-Mitra} to groups
with left invariant metrics. Section \ref{SuperNilpotent} is devoted to an intermediate
step in our proof: the study of the Assouad-Nagata dimension of nilpotent groups. In
Section \ref{MainSection} we will prove that the Assouad-Nagata dimension of
simply-connected solvable Lie groups is equal to the topological dimension. As a consequence we will compute the asymptotic Assouad-Nagata dimension of finitely generated polycyclic groups. In section
\ref{MainSection2} we will study, using Iwasawa decompositions, the Assouad-Nagata
dimension of semisimple Lie groups. Finally in section \ref{MainSection3} we will prove
the main result of the paper. \\

{\bf Acknowledgements:} Both authors are very grateful to the Hausdorff Research
Institute for Mathematics and specially with the organizers of the 'Rigidity program' for
their hospitality and support. They also thank David Fisher for useful conversations and
U. Lang for very helpful comments and suggestions, specially for pointing out \cite{Car}
and the necessity of the example \ref{ExampleCarlsson}.

\section{Assouad-Nagata dimension and transformations on the metric}\label{Nagata}
Let $s$ be a positive real number. A {\it $s$-scale chain} (or $s$-path)
 between two points $x$ and $y$ of a metric space $(X, d_X)$ is defined as
a finite sequence of points
$\{x= x_0, x_1, ..., x_m = y\}$ such that
$d_X(x_i, x_{i+1}) < s$ for every $i = 0, ..., m-1$.
A subset $S$ of a metric space $(X, d_X)$ is said to be {\it $s$-scale connected} if every two elements of $S$ can be connected by $s$-scale chain
contained in $S$.
\begin{Def} A metric space $(X, d_X)$ is said to be of
{\it asymptotic dimension} at most $n$ (notation $\asdim(X, d) \le
n$) if there is an increasing function $D_X: \RR_+ \to \RR_+$ such
that for any $s> 0$ there is a cover $\UU =\{\UU_0, ...,\UU_n\}$
so that the $s$-scale connected components of each $\UU_i$ are
$D_X(s)$-bounded i.e. the diameter of every component is bounded
by $D_X(s)$. \end{Def}
 The function $D_X$ is called an {\it
$n$-dimensional control function} for $X$. Depending on the type
of $D_X$ one can define the following invariants:

\begin{Def} A metric space $(X, d_X)$ is said to have
\begin{itemize}
\item {\it  Assouad-Nagata dimension} at most $n$ (denoted by $\dim_{AN} (X, d) \le n$) if it
has an $n$-dimensional control function $D_X$ of the form $$D_X(s) = C\cdot s$$ with $C>0$ some fixed constant.

\item {\it asymptotic Assouad-Nagata dimension} at most $n$ (denoted by
    \newline $\asdim_{AN} (X, d) \le n$) if it has an $n$-dimensional control function
    $D_X$ of the form
    $$D_X(s) = C\cdot s +k$$ with $C >0$ and $k\in \RR$ two fixed constants.

\item {\it capacity dimension} at most $n$ (notation $\cdim(X, d) \le n$) if it has an $n$-dimensional
control function $D_X$ such that $$D_X(s) = C\cdot s$$ in a neighborhood of $0$ i.e. with $s$ sufficiently small. \end{itemize} \end{Def}

\begin{Rem}
 For any metric space $(X, d)$ we have $\dim(X, d) \le dim_{AN}(X, d)$ where $\dim$ is the topological dimension (see \cite{Assouad}).
\end{Rem}
One map $f: (X, d) \to (Y, D)$ between metric spaces is said to be a {\it quasi-isometric embedding} if there exists two constants $C \ge 1, \lambda \ge 0$ such that:
\begin{equation}\label{quasiisometry}
\frac{1}{C}\cdot d(x, y) - \lambda \le D(f(x), f(y)) \le C\cdot d(x, y) + \lambda.
\end{equation}
If in addiction there exists a $K >0$ such that $D(y, f(X)) \le K$ for every $y \in Y$,
$f$ is said to be a {\it quasi-isometry} and the spaces $(X, d)$ and $(Y, D)$ are said to
be quasi-isometric. In (\ref{quasiisometry}), if $\lambda = 0$ then $f$ is said to be a {\it bilipschitz equivalence} and the
spaces are said to be {\it bilipschitz equivalent}. One important fact about Assouad-Nagata dimension is that it is
preserved under bi-Lipschitz equivalences. Analogously the asymptotic Assouad-Nagata
dimension is preserved under quasi-isometries. The capacity dimension is invariant under
maps that are bi-Lipschitz at small scales i.e. bi-Lipschitz when the distances are less
than some fix constant $\epsilon>0$. \medskip

It is clear that $\asdim(X, d)  \le \asdim_{AN}(X, d)$. The relation
among the Assouad-Nagata dimension, the asymptotic Assouad-Nagata
dimension and the capacity dimension in a metric space was studied
in \cite{Brod-Dydak-Higes-Mitra} via two functors $\max(d, 1)$ and
$\min(d, 1)$, one in the category of bounded metric spaces and the
other in the category of discrete metric spaces. We collect from
that paper the most important results for our purposes in the
following.

\begin{Prop}\label{MicroMacro}\cite{Brod-Dydak-Higes-Mitra}
Let $(X, d)$ be a metric space. Then:
\begin{enumerate}
\item $\cdim(X, d) = \dim_{AN}(X, \max(d, \epsilon))$ for every $\epsilon >0$.
\item $\asdim_{AN}(X, d) = \dim_{AN}(X, \min(d, \epsilon))$ for every $\epsilon >0$.
\item $\dim_{AN} (X,d) =
max\{\asdim_{AN} (X, d), \cdim (X,d)\}$
\end{enumerate}
\end{Prop}

It was shown in $\cite{Lang}$ that the Assouad-Nagata dimension is also invariant under
quasisymmetric embeddings as for example snow-flake transformations. Recall that a {\it
snow flake transformation} of a metric space $(X, d)$ is of the form $(X , d^{\alpha}) $
where $0 < \alpha \le 1$. We will need the following lemma that is in some sense a
generalization of this fact for geodesic spaces:

\begin{Lem}\label{Traductor2}Let $(X, d)$ be a (non bounded) geodesic metric space. Suppose $f: \RR_+ \to \RR_+$ is a surjectively non decreasing function such that $(X, f(d))$ is a metric space. Then $\dim_{AN}(X, f(d)) \le \dim_{AN}(X, d)$.
\end{Lem}
\begin{pf}
Notice that $f(x+y) \le f(x)+f(y)$ for every $x, y \in \RR_+$. This follows from the fact
that $(X, d)$ is geodesic, unbounded, and $(X, f(d))$ is a  metric. Let us show it. Pick
$a, b \in X$ such that $d(a, b) = x+y$. As $(X, d)$ is geodesic there exists a $c \in X$
such that $d(a, c) = x$ and $d(b, c) = y$. Hence by the triangle inequality in $(X,
f(d))$ we get $f(d(a, b)) \le f(d(a, c)) + f(d(b, c))$. \medskip

Suppose $dim_{AN}(X, d) \le n$. Let $D(s) = C \cdot s$ be an $n$-dimensional control
function of $(X, d)$. Without loss of generality we can assume $C \in \NN$. Let $s>0$.
Take some inverse $f^{-1}$ of $f$. Hence there exists a covering of $(X, d)$ of the form
$\UU = \bigcup_{i= 0}^n \UU_i$ such that the $f^{-1}(s)$-scale components of each $\UU_i$
are $C \cdot f^{-1}(s)$-bounded for every $i \in \{0, \cdots, n\}$. This implies that if
$x, y \in X$ are in different $f^{-1}(s)$-scale components of $\UU_i$ then $s \le f(d(x,
y))$. Therefore the $s$-scale components of $\UU_i$ in $(X, f(d))$ are contained in the
$f^{-1}(s)$-scale components of $\UU_i$ in $(X, d)$. Let $x, y \in \UU_i$ be two elements
that belongs to the same $f^{-1}(s)$-scale component of $\UU_i$. We have $d(x, y) \le
C\cdot f^{-1}(s)$. This implies $f(d(x, y)) \le f(C \cdot f^{-1}(s))$. On the other hand
by the subadditivity of $f$  and $C \in \NN$ we get $ f(C\cdot f^{-1}(s) \le C\cdot
f(f^{-1}(s)) = C \cdot s$ as desired. \end{pf}

From the previous lemma we have the following result that can be
applied to many remarkable cases as for example unbounded trees or
Cayley graphs of finitely generated groups such that the asymptotic
dimension coincides with the asymptotic Assouad-Nagata dimension.

\begin{Cor} Let $(X, d)$ be a geodesic space and let $f: \RR_+ \to \RR_+$ as in the previous lemma such that  $lim_{x \to \infty} f(x) = \infty$. If $\asdim(X, d) = \dim_{AN}(X, d)$ then $\dim_{AN}(X, f(d)) = \dim_{AN}(X, d)$
\end{Cor}
\begin{pf}
Just notice that $\asdim(X, d) = \asdim(X, f(d))$ as $f$ induces a coarse equivalence.
\end{pf}
The condition on the equality of dimensions can not be dropped as it is shown in the
following example.

\begin{ex} In \cite{Brod-Dydak-Lang} it was shown that the Cayley graph $\Gamma$ of the group $\ZZ_2 \wr \ZZ^2$
has infinite Assouad-Nagata dimension but asymptotic dimension
equals two. Moreover from Theorem 5.5 of \cite{Brod-Dydak-Lang} it
follows that the Cayley graph has a $2$-dimensional control function
of polynomial type. Let $d$ be the metric of the Cayley graph
$\Gamma$. It is not hard to check (or you can apply directly
proposition 2.2 of \cite{Higes}) that $dim_{AN}(\Gamma, log(d+1)) =
2$.
\end{ex}

\section{Dimension of exact sequences}\label{Hurewicz}

In this section we extend  a result of \cite{Brod-Dydak-Levin-Mitra}
about dimension of exact sequences of finitely generated groups
equipped with word metrics to dimension of exact sequences of
general groups(not necessarily countable) equipped with left
invariant metrics. First we will recall some definitions and results
from \cite{Brod-Dydak-Levin-Mitra}.

\begin{Def}
Given a function between metric spaces $f: X \to Y$ and given $m \in
\NN$. An $m$-dimensional control function of $f$ is a function $D_f:
\RR_+ \times \RR_+ \to \RR_+$ such that for all $r_X >0$ and $R_Y
>0$, every  set  $A \subset X$ such that $\diam(f(A)) \le R_Y$ can be
expressed as the union of $m+1$ subsets $\UU = \bigcup_{i= 0}^n
\UU_i$ such that the $r_x$-scale components of each $\UU_i$ are
$D_f(r_X, R_Y)$-bounded.
\end{Def}
As in the previous section depending on the type of $D_f$ we could
have different notion of dimension of functions (see
\cite{Brod-Dydak-Levin-Mitra}). We will use just the following:

\begin{Def} A function
between metric spaces $f: X \to Y$ is said to have Assouad-Nagata
dimension at most $n$ (notation $\dim_{AN}(f) \le n$) if there
exists an $n$-dimensional control function $D_f$ of the form
$D_f(r_x, R_Y) = a \cdot r_x + b \cdot R_Y$. If $n$ is the minimum
number such that $f$ satisfies this property then $f$ is said to
have Assouad-Nagata dimension exactly $n$.
\end{Def}

\begin{Prop}\cite[Theorem 7.2]{Brod-Dydak-Levin-Mitra}\label{ANfunction2}
If $f: X \to Y$  is a Lipschitz function between metric spaces then:
\begin{equation}\label{ANfunction}
\dim_{AN}(X, d_X)\le \dim_{AN}(f) + \dim_{AN}(Y, d_Y).
\end{equation}
\end{Prop}
\begin{Prop}\cite[Proposition 3.7]{Brod-Dydak-Levin-Mitra}\label{PropSubset}
Suppose $A$ is a subset of a metric space $(X, d)$, $m \ge 0$, $R >0$. If $D_A$ is an
$m$-dimensional control function of $A$ then $D_B(s) := D_A(s +2R) +2R$ is an
$m$-dimensional control function of the $R$-neighborhood $B(A, R)$. \end{Prop}

Let $1 \to K \to G \to H \to 1$ be an exact sequence of groups. Suppose $d_G$ is a left
invariant metric. This metric induces in $H$ a natural metric called the {\it Hausdorff
metric} defined by the norm $\|g\cdot K\|_H = inf\{\|g\cdot k\|_G; k \in K\}$. Notice
that when $G$ is a finitely generated group and $d_G$ is a word metric then the
corresponding induced metric $d_H$ is also a word metric.

The induce Hausdorff metric allows us to extend Corollary 8.5 of
\cite{Brod-Dydak-Levin-Mitra} to general groups as follows.

\begin{Prop}\label{ExactSequence}
Let $1 \to K \to G \to H \to 1$ be an exact sequence of groups.  Then
\begin{enumerate}
\item $\dim_{AN}(G, d_G) \le \dim_{AN}(K, d_G|_K) + \dim_{AN}(H, d_H);$
\item $\asdim_{AN}(G, d_G) \le \asdim_{AN}(K, d_G|_K) + \asdim_{AN}(H, d_H);$
\item $cdim(G, d_G) \le cdim(K, d_G|_K) + cdim(H, d_H),$ \end{enumerate} where $d_G$
    is a left invariant metric of $G$ and $d_H$ is the induced (left invariant)
    Hausdorff metric on $H$.
\end{Prop}
\begin{pf}
We will prove only the first inequality. The other two follow from the first one and
Proposition \ref{MicroMacro} (notice that the induced Hausdorff metrics are preserved by
the two functors). \smallskip

The proof is close to that of Corollary 8.5 in \cite{Brod-Dydak-Levin-Mitra}.  First, it
is clear that if $d_H$ is the Hausdorff metric, then the projection map $f: (G, d_G) \to
(H, d_H)$ is $1$-Lipschitz. Let $B = B(1_H, R_H)$ and let $A = f^{-1}(B)$.  Take $a \in
A$, then $\|f(a)\|_H < R_H$. As $f$ is an homomorphism and by the symmetry of the norm we
get $\|f(a^{-1}\|_H = \|f(a)^{-1}\|_H = \|f(a)\|_H < R_H$. But on the other hand we have
$\|f(a^{-1})\|_H = inf\{\|a^-1\cdot k\|_G; k \in K\}$. Therefore $\|f(a^{-1})\|_H =d_G(a,
K) < R_H$, and we have shown $A \subset B_G(K, R_H)$.  By Proposition \ref{PropSubset} we
have that $D_f(r_G, r_H) := D_K(r_G + 2\cdot R_H) + 2 \cdot R_H$ is an $m$-dimensional
control function of $B(K, R_H)$, provided $D_K$ is an $m$-dimensional control function of
$K$.  Now if $dim_{AN}(K, d_G|_K) \le m$ then there exists an $C >0 $ such that $D_K (s)
= C \cdot s$. Since all $f^{-1}(B(y, R_H))$ are isometric we can apply Proposition
\ref{ANfunction2} to get the inequality:  $$dim_{AN}(G, d_G)\le dim_{AN}(K, d_G|_K) +
dim_{AN}(H, d_H).$$  \end{pf}

\section{Assouad-Nagata dimension of Nilpotent Lie groups}\label{SuperNilpotent}

As an intermediate step we study in this section the Assouad-Nagata
dimension of Nilpotent Lie groups. From Proposition \ref{MicroMacro} we would need to calculate the capacity
dimension and the asymptotic Assouad-Nagata dimension. For all
connected Lie groups (not necessarily solvable) we will show first that
the capacity dimension is less or equal than the topological one. We
will use the ideas of Buyalo and Lebedeva. Following definitions and
theorem \ref{BuyaloThm} come from \cite{Buyalo}:
\begin{Def}
A map $f:Z\to Y$ between metric spaces is said to be {\em
quasi-homothetic} with coefficient $R>0$, if for some $\lambda\ge 1$
and for all $z$, $z'\in Z$, we have
$$R\cdot d_Z(z, z')/\lambda \le d_Y(f(z),f(z')) \le \lambda \cdot R\cdot d_Z(z, z').$$
In this case, it is also said that $f$ is $\la$-quasi-homothetic with
coefficient $R$.
\end{Def}

\begin{Def}
A metric space $Z$ is {\em locally similar} to a metric space $Y$,
if there is $\la\ge 1$ such that for every sufficiently large $R>1$
and every $A\subset Z$ with $\diam A\le\La_0/R$, where
$\La_0=\min\{1,\diam Y/\la\}$, there is a $\la$-quasi-homothetic map
$f:A\to Y$ with coefficient $R$.
\end{Def}

\begin{Thm}\cite[Theorem 1.1]{Buyalo} \label{BuyaloThm} Assume that a metric space
$Z$ is locally similar to a compact metric space $Y$. Then $\cdim
(Z)<\infty$ and $\cdim (Z)\le\dim (Y)$.
\end{Thm}

\begin{Lem}
Let $G$ be a connected Lie group of dimension $n$ with a
left-invariant Riemannian metric $d_{G}$.  Then it is locally
similar to a closed ball in $\RR^{n}$.\end{Lem}

\begin{pf}
Note that on a smooth manifold $M$, the exponential map $Exp: T_{g}M
\rightarrow M$ is a diffeomorphism between a closed ball of radius
$\epsilon$ (which might depend on the point $p$) and a convex
neighborhood of $p$. Now, any diffeomorphism between two compact
Riemannian metric spaces is bilipschitz, so there exists a constant
$\lambda$ (which might also depend on the point $p$) such that the
exponential map restricted to the ball of radius $\epsilon$ is
$\lambda$-bilipschitz for some $\lambda \ge 1$. If $M$ is a Lie
group, by left translations neither the $\epsilon$  nor the
$\lambda$ depends on the point $p$, so we can assume $p = 1_G$. Let
$\log: \bar{B}(1_{G}, \epsilon) \rightarrow \bar{B}(0,K)$ be the
inverse of the exponential map where  $\bar{B}(0,K)$ is some ball of
$\RR^n$ that contains the image of $\log(\bar{B}(1_{G}, \epsilon))$.
The metric $d$ of $\RR^n$ that we are considering is the usual one.
Without loss of generality let us suppose $K
> 2\cdot \lambda$. In such case $\Lambda_0 := min\{1, \frac{2\cdot
K}{\lambda}\} = 1$ and $R_0 = \frac{1}{\epsilon}$. Given $R > R_0$,
let $A$ be a closed ball of radius at most $\frac{1}{R}$. By left
translation we can assume $A$ is centered in $1_G$. Notice
$\log(1_G) = 0$. Define the natural dilatation $g: \bar{B}(0, K) \to
\bar{B}(0, K \cdot  R)$. By the assumption of $K > 2 \cdot \lambda$
the image of the composition of $g\circ \log$ restricted to $A$ lies
in $\bar{B}(0, K)$. Moreover it satisfies:
$$\lambda^{-1}\cdot d_G(x, y) \le 1/R\cdot d(g(\log(x)), g(\log(y))) \le
\lambda \cdot d_G(x, y).$$
\end{pf}

Combining this lemma with Theorem \ref{BuyaloThm} we get:

\begin{Cor}\label{CapDimLie}
Let $G$ be a connected Lie group equipped with a left-invariant Riemannian metric
$d_{G}$. Then $cdim(G, d_G) \le dim(G)$.\end{Cor}

Now we study the asymptotic Assouad-Nagata dimension of a nilpotent Lie group $G$ . For
such purpose we will apply Proposition \ref{ExactSequence} to the exact sequence $1 \to
N^r \to G \to G/N^r \to 1$ where $N^r$ is the last term in the lower central series.
Notice that $G/N^r$ is also a nilpotent Lie group.  The metrics considered will be a
Riemannian metric in $G$, $d_G$, and the corresponding induced Hausdorff metric in $G/N^r$.
We will need the following.

\begin{Lem}\label{RieInQuotient}
Let $G$ be a connected Lie group and $H \lhd G$ a normal subgroup such that $G/H$ is a
connected Lie group.  Then there are Riemannian metrics $\rho_{G}$ and $\rho_{G/H}$ on
$G$ and $G/H$ such that the (left-invariant) Hausdorff metric on $G/H$ induced by
$\rho_{G}$ agrees with the path metric induced by $\rho_{G/H}$. \end{Lem}

\begin{pf}
Write $\pi: G \rightarrow G/H$ for the canonical projection.  In $G$, choose a complement
$V$ to $T_{e}H$ in $T_{e}G$.  Define $\rho_{G}$ at $T_{e}G$ by choosing inner products on
$T_{e}H$ and $V$, and define $\rho_{G/H}$ at $T_{e}(G/H)$ to be the inner product chosen
on $V$.  By left translating $V$ we obtain a distribution $\triangle_{V}$ on $G$, so
every differentiable path $\gamma \in G/H$ has a isometric lift $\tilde{\gamma}$ lying
entirely in $\triangle_{V}$, which is unique up to the choice of the starting point.
\smallskip

Fix a coset $pH \in G/H$.  Since both metrics are left invariant, it suffices to show
that $d_{\rho_{G/H}}(H, pH)$ is the same as $d_{\mathcal{H}}(H, pH)$.  Let $\mathcal{K}$
be the set of isometric lifts of differentiable paths connecting $H$ with $Hp$ in $G/H$
that start at the identity, and $\mathcal{B}$ be the set of differentiable paths starting
at the identity and end at some point in $pH$.  That is, \[ {\mathcal{K}}=\{
\tilde{\gamma}: \tilde{\gamma}(0)=e, \pi(\tilde{\gamma})=\gamma, \gamma:[0,L] \rightarrow
G/H, \gamma(0)=H, \gamma(L)=Hp \}\] Clearly $\mathcal{K} \subset \mathcal{B}$. However if
$\eta \in \mathcal{B}$, $\widetilde{\pi(\eta)}$, the lift of its projection, is also in
$\mathcal{B}$, with length no bigger than $\eta$, so $\min \{ \| \eta \|: \eta \in
\mathcal{B} \}=\min \{ \| \widetilde{\pi(\eta)} \|: \eta \in \mathcal{B} \}.$  But $$ \{
\widetilde{\pi(\eta)}: \eta \in \mathcal{B} \} = \mathcal{K}.$$  The claim now follows
since \[ d_{{\mathcal{H}}}(H,pH) =\min \{ \| \eta \|: \eta \in {\mathcal{B}} \} = \min \{ \|
\zeta \|: \zeta \in {\mathcal{K}} \} =d_{\rho_{G/H}}(H, Hp) \] \end{pf}

Let $N$ be a connected, simply connected nilpotent Lie group, and
let $N^{i}$ be the $i$-th term in its lower central series. That is,
$N^{2}=[N,N]$, and $N^{i+1}=[N, N^{i}]$. By construction, each
quotient $N^{i}/N^{i+1}$ is an abelian Lie group of dimension $n_{i}$, so by fixing a
subset ${\mathcal{K}}_{i} \subset N^{i}$, a set of $N^{i+1}$ coset
representatives in $N^{i}$, we have a bijection $\phi_{i}:
{\mathbb{R}}^{n_{i}} \rightarrow {\mathcal{K}}_{i}$, and a map
$\phi: N \rightarrow \oplus_{i} {\mathbb{R}}^{n_{i}}$ defined as
\[ p \stackrel{\phi}{\longrightarrow} (p_{1}, p_{2} \cdots ) , p_{i} \in {\mathbb{R}}^{n_{i}} \] \noindent where
$\phi_{1}(p_{1})N^{2}=pN^{2}$, and $\phi_{i}(p_{i}) N^{i+1} = (\phi_{i-1}(p_{i-1}))^{-1}
\cdots (\phi_{2}(p_{2}))^{-1} (\phi_{1}(p_{1}))^{-1} p N^{i+1}$.  One can check that
$\phi$ is a bijection with the inverse given by
\[ (p_{1}, p_{2}, \cdots ) \stackrel{\phi^{-1}}{\longrightarrow} \phi_{1}(p_{1})
\phi_{2}(p_{2}) \phi_{3}(p_{3}) \cdots \]

\noindent Note that if the degree of nilpotency of $N$ is $r$, then $N^{r}$, the last
non-trivial satisfies $$\phi(N^{r}) =  \{ (0,0, \cdots, p_{r}): p_{r} \in
{\mathbb{R}}^{n_{r}} \}.$$

With this coordinate system we define $D: N \times N \rightarrow {\mathbb{R}}$ as $$
D(p,q)=D(1,p^{-1}q) = \sum_{i}^{r} \| \vec{x}_{i} \|^{1/i}, $$ where
$\phi(p^{-1}q)=(\vec{x}_{1}, \cdots \vec{x}_{r})$, and $\|  \|$ is the
standard Euclidean norm on ${\mathbb{R}}^{n}$. \\

The main result from \cite{Karidi} is the following.

\begin{Thm}\cite[Theorem 4.2]{Karidi} \label{Karidi} Let $N$ be a
connected, simply connected nilpotent Lie group.  Then there is a Riemannian metric
$d_{N}$ on $N$ and constant $\kappa $ such that for any point $p$ with $d_{N}(e, p) > 1$,
$$ 1/\kappa D(1, p) \leq d_{N}(1, p) \leq \kappa  D(1,p).$$ Equivalently, $(N, \min(1,
D))$ is bilipschitz to $(N, \min(1, d_{N}))$. \end{Thm}

An understanding of distance distortions in nilpotent Lie groups can now be obtained.

\begin{Lem} \label{DistorsionH}
Let $N$ be a connected, simply connected nilpotent Lie group, and $H$ the last
non-trivial term in its lower central series.  Then there are Riemannian metrics $d_{N}$,
$d_{H}$ on $N$ and $H$ such that $(H, \min(1,d_{N}|_{H}))$ is bilipschitz to $(H,
\min(1,(d_{H})^{1/r}))$ where $r$ is the degree of nilpotency of $N$.  \end{Lem}

\begin{pf}
Note that $H=N^{r}$.  Equip $H$ with the Riemannian metric induced by $\phi_{r}:
{\mathbb{R}}^{n_{r}} \rightarrow H$, while putting on $G$ the Riemannian metric given by
Theorem \ref{Karidi}. With these choices we have $D|_{H}= (d_{H})^{1/r}$, and the claim
follows since Theorem \ref{Karidi} says that $(H, \min(1, D|_{H}) )$ is bilipschitz to
$(H, \min(1, d_{N}|_{H}))$.  \end{pf}

\begin{Lem}\label{NilpotentCase}
Let $N$ be a nilpotent Lie group and let $d_N$ be a left invariant Riemannian metric.
Then $asdim_{AN}(N, d_N) \le dim(N)$ \end{Lem}

\begin{pf}
Recall the topological dimension of $N$ is the same as the sum of the topological
dimensions of factors in its lower central series. We will prove the lemma by induction
on the degree of nilpotency. The base case is when $N$ is an abelian Lie group, and in
this case we have the equality $asdim_{AN}(N, d_{N})=dim(N)$.  In general we consider the
following short exact sequence $$ 1 \to H \to N \to N/H \to 1 $$where $H$ is the last
term in the lower central series of $N$. Let $d_{N/H}$ denotes the induced Hausdorff
metric on $N/H$ from $d_N$ and $d_{H}$ denotes the Riemannian metric on $H$ induced by
$d_{N}$.  By Proposition \ref{ExactSequence} we have: $$\asdim_{AN}(N, d_N) \le
\asdim_{AN}(H, d_N|_H) + \asdim_{AN}(N/H, d_{N/H}).$$

But $\asdim_{AN}(H, d_N|_H) = \asdim_{AN}(H, d_H)$ because of Lemma \ref{DistorsionH} and
the fact that the Assouad-Nagata dimension is invariant under quasi-isometries and
snowflake transformations.  Moreover $\asdim_{AN}(H, d_{H})=dim(H)$ since $H$ is abelian.
 On the other hand, Lemma \ref{RieInQuotient} says that $d_{N/H}$ is a Riemannian metric
in $N/H$, and since $N/H$ is nilpotent with one less degree of nilpotency, induction
hypothesis yields $asdim_{AN}(N/h, d_{N/H}) \le dim(N/H)$, to which the desired claim now
follows. \end{pf}

\begin{Thm}\label{NilpotentCase2}
Let $N$ be a nilpotent Lie group and let $d_G$ be some Riemannian metric.  Then
$dim_{AN}(N, d_N) = dim(N)$ \end{Thm}

\begin{pf}
By Proposition \ref{MicroMacro} the Assouad-Nagata dimension is equal the maximum of the capacity dimension and the asymptotic Assouad-Nagata dimension. Hence applying  corollary
\ref{CapDimLie} and lemma \ref{NilpotentCase} we get $dim_{AN}(N,
d_N)\le dim(N)$. The other inequality follows from the fact that the Assouad-Nagata dimension is
always greater or equal to the topological one (see \cite{Assouad}).
\end{pf}

Next example will show that the proof of Theorem 3.5 of \cite{Car} can not be applied to
compute the Assouad-Nagata dimension of nilpotent groups. We recommend the reading of
such proof in order to understand better the example.

\begin{ex}\label{ExampleCarlsson}
Consider the 4-dimensional nilpotent Lie group determined by the following Lie algebra.
\[ [e_{1},e_{2}]=e_{3}, [e_{1}, e_{3}]=e_{4} \]

\noindent It has a group structure given by
\[ \left(
     \begin{array}{c}
       x_{1} \\
       x_{2} \\
       x_{3} \\
       x_{4} \\
     \end{array}
   \right) \bullet \left(
                     \begin{array}{c}
                       y_{1} \\
                       y_{2} \\
                       y_{3} \\
                       y_{4} \\
                     \end{array}
                   \right) = \left(
                               \begin{array}{c}
                                 x_{1}+y_{1} \\
                                 x_{2}+y_{2} \\
                                 x_{3}+y_{3} + \frac{1}{2}(-x_{2}y_{1}+x_{1}y_{2})  \\
                                 x_{4}+y_{4} + \frac{1}{12}(x_{1}-y_{1})(-x_{2}y_{1}+x_{1}y_{2}) + \frac{1}{2}(-x_{3}y_{1}+x_{1}y_{3}) \\
                               \end{array}
                             \right) \]

A left-invariant Riemannian metric is given by
\[ \left( dx_{1} + \frac{1}{2} x_{2} dx_{3} -\frac{1}{6} x_{1}^{2} x_{2} dx_{4}
\right)^{2} + \left( dx_{2} - \frac{1}{2} x_{1} dx_{3} + \frac{1}{6} x_{1}^{2} dx_{4} \right)^{2}
+ \left( \frac{-1}{2} x_{1} dx_{3} \right)^{2} + dx_{4}^{2} \]

Following the proof of Theorem 3.5 in \cite{Car}, we express this 4-dimensional nilpotent
groups as a semidirect product $T \ltimes N_{0}$, where $N_{0}$ is the subgroup generated
by $e_{2},e_{3},e_{4}$, and $T$ is the subgroup generated by $e_{1}$.  Then the metric on
$N_{0}(x_{1}e_{1})$ is

\[ \left( dx_{2}  + x_{1}dx_{3} + x_{1}^{2} dx_{4} \right)^{2} + \left( dx_{3} +
x_{1}dx_{4} \right)^{2} + \left( dx_{4} \right)^{2} \]

\noindent Note the role of $e_{1}$ coordinate, $x_{1}$, plays in the metric.  For example
if $I$ is an interval of size $c$ in the $e_{4}$ direction, then the diameter of
$I(x_{1}e_{1})$, right translate of $I$ by $x_{1}e_{1} \in T$, is $x_{1}^{2}c$. Similarly
if $J$ is an interval of size $c$ in the $e_{3}$ direction then the diameter of
$J(x_{1}e_{1})$ is $x_{1}c$. In this way we see for a subset $W \subset N_{0}$ of bounded
diameter, the diameter $W(x_{1}e_{1})$ depends on the value of $x_{1}$.  Since the metric
is left-invariant, the diameter of $W(x_{1}e_{1})$ is the same as the diameter of
$(x'_{1}e_{1})W(x_{1}e_{1})$ for any $x'_{1}e_{1} \in T$, the same is true for conjugates
of $W$ by $x_{1}e_{1}$. \smallskip

Now to see that the diameter of $S_{i}^{l}(U)$ is about $D^{2}$ times the diameter of
$U$, we observe that $S_{i}^{l}(U)$ consists of right translates of $U$ by elements of
$T$.  Express $U$ as $U=t_{i}W$ for some $W \subset N_{0}$, $t_{i} \in T$, we have that
$$S_{i}^{l}(U)= \bigcup_{t \in 2D} Ut =\bigcup_{t \in 2D}  t_{i}Wt =\bigcup_{t \in 2D}
t_{i}t (t^{-1}Wt) ,$$ which makes the diameter of $S_{i}^{l}(U)$ the maximum of $2D$ and
the diameter of $t^{-1}Wt$.  For $D<1$, the diameter of $S_{i}^{l}(U)$ varies linearly
with $D$ (since $D> D^{2}$ for $D< 1$), but for $D > 1$, the discussion above shows that
the diameter depends on $D^{2}$.  \end{ex}

\section{Assouad-Nagata dimension of solvable Lie groups}\label{MainSection}
\begin{Def}
Let $G$ be a connected, simply connected solvable Lie group. The exponential radical,
denoted as $Exp(G)$, is a closed normal subgroup such that $G/Exp(G)$ is the biggest
quotient with polynomial growth. \end{Def}

Osin shown in \cite{Osin} that given a Riemannian metric $d_G$ on  $G$ and a Riemannian
metric on the exponential radical $d_{Exp(G)}$ then there exists two constants $C>0$ and
$\epsilon \ge 0$ such that fore every $h\in Exp(G)$:
\begin{equation}\label{OsinInequality}
\frac{1}{C}log(\|h\|_{Exp(G)}+1) - \epsilon \le \|h\|_G\le C\cdot
log(\|h\|_{Exp(G)}+1)+\epsilon
\end{equation}

Also he proved that the exponential radical is contained in the nilradical of $G$, so it
is nilpotent and we have $asdim_{AN}(Exp(G), d_{Exp(G)}) \le dim(Exp(G)))$ by Lemma
\ref{NilpotentCase}.  But if we apply Lemma \ref{Traductor2} to the particular case $f(s)
= log(s+1)$ and the invariance of the asymptotic Assouad-Nagata dimension under
quasi-isometries we get $asdim_{AN}(Exp(G), d_G|_{Exp(G)}) \le dim(Exp(G))$. On the other
hand by the definition of exponential radical there is an exact sequence:
$$1 \to Exp(G) \to G \to S \to 1,$$ where $S$ is a solvable Lie
group with polynomial growth, and quasi-isometric to a nilpotent group of the same
topological dimension \cite{Breuillard}.  Hence by Lemma \ref{NilpotentCase} and Lemma
\ref{RieInQuotient} we have $\asdim_{AN}(S, d_S) \le dim(S)$ where $d_S$ is the
Riemannian metric induced by $d_G$. Therefore applying Proposition \ref{ExactSequence} we
get:

\begin{Prop}\label{SolvAN}
Let  $G$ be a connected solvable Lie group then:
$$\asdim_{AN}(G, d_G) \le \asdim_{AN}(Exp(G), d_1) + \asdim_{AN}(G/Exp(G), d_2) \le dim(G)$$
where $d_1$ and $d_2$ are two Riemannian metrics defined in $Exp(G)$
and $G/Exp(G)$ respectively. \end{Prop}

Now we can get the main result of this paper:

\begin{Thm}\label{MainTheorem} Let $G$ be a connected solvable Lie group then:
$$\dim_{AN}(G, d_G) = dim(G).$$ \end{Thm}
\begin{pf}
The proof is analogous to the that of Theorem \ref{NilpotentCase2}.  On the one hand we
have $dim_{AN}(G, d_G) \ge dim(G)$ by \cite{Assouad}.  By Proposition \ref{MicroMacro}
the Assouad-Nagata dimension is equal the maximum of the capacity dimension and the
asymptotic Assouad-Nagata dimension. Hence the other inequality follows from Proposition
\ref{SolvAN} and Corollary \ref{CapDimLie}.
\end{pf}

As a consequence of this theorem and  Theorem 1.3 of \cite{Lang} we
get the following interesting property of connected solvable Lie
groups:
\begin{Cor}
Let $G$ be a connected solvable Lie group equipped with a Riemannian
metric $d_{G}$.  Then there exists a $0 < p \le 1$ such that $(G,
d_{G}^{p})$ can be bilipschitz embedded into a product of $dim(G)+1$
many trees. \end{Cor}

\begin{Def}
Let $\Gamma$ be a finitely generated solvable group.  Then the
hirsch length is defined as:
$$h(\Gamma) = \sum \dim_{\QQ}( \Gamma_{i}/\Gamma_{i+1} \otimes \QQ)$$
where $\Gamma_{2}=[\Gamma, \Gamma]$ and $\Gamma_{i+1} =
[\Gamma_{i},\Gamma_{i}]$ \end{Def}

Next Corollary answers Question 4 of \cite{Nowak} and problem 8.3 of \cite{DranProblems}:

\begin{Cor}
Let $(\Gamma, d_w)$ be a polycyclic group equipped with a word metric $d_{w}$.  Then
$$\asdim_{AN}(\Gamma, d_w) = h(\Gamma)$$ \end{Cor}
\begin{pf}
In \cite{Dran-Smith} it was proved that $\asdim(\Gamma, d_w) \geq h(\Gamma)$, so the
unique thing we need to show is that $\asdim_{AN}(\Gamma, d_w) \le h(\Gamma)$.  It is
known that every polycyclic group is a cocompact lattice in a connected, simply connected
solvable Lie group $H$ after modding out a finite torsion subgroup. Moreover $h(\Gamma) =
\dim(H)$ and by Theorem \ref{MainTheorem} we have $\dim(H) = \dim_{AN}(H)$. As $\Gamma$
is a lattice we have $\asdim_{AN}(\Gamma) = \asdim_{AN}(H)\le \dim_{AN}(H)$, where the
last inequality follows from Proposition \ref{MicroMacro}.  \end{pf}

\begin{Rem}
The condition of word metrics can not be relaxed as it is shown in \cite{Higes2}, that a
countable nilpotent group $G$ can always be equipped with a proper left invariant metric
$d_G$,  such that $\asdim_{AN}(G, d_G)$ is infinite. \end{Rem}

\section{Assouad-Nagata dimension of semisimple Lie groups}\label{MainSection2}
In this section we will compute the asymptotic Assouad-Nagata dimension of semisimple Lie
groups. The idea will be to study the Iwasawa decompositions of such groups. Theorem
\ref{Iwasawa} from \cite{Knapp} and \cite{Hochschild} gives the structural description of
a semisimple Lie group that we will need.

\begin{Def}
A Lie algebra is semisimple if it does no have a non-trivial solvable ideal. \end{Def}

\begin{Def}
A Lie group is semisimple if its Lie algebra is semisimple. \end{Def}

\begin{Thm}\textit{Theorem 6.31, 6.46 from \cite{Knapp}, Theorem 3.1, Lemma 3.3, Chap XV from
\cite{Hochschild}}\label{Iwasawa} \newline

Let $G$ be a semisimple Lie group with finitely many component.  Then there exist
subgroups $K$, $A$ and $N$ such that the multiplication map $A \times N \times K
\rightarrow G$ given by $(a,n,k) \mapsto ank$ is a diffeomorphism. The groups $A$ and $N$
are simply connected abelian and simply connected nilpotent respectively, and $A$
normalizes $N$.\smallskip

Furthermore, $Z(G)$, the center of $G$, is contained in $K$, and there is an
isometrically embedded connected, simply connected abelian Lie group $V = Z(K)$, a
compact subgroup $T < K$ such that $K$ is the semidirect product $V \rtimes T$, and that
$Z(G)/V$ is compact. \smallskip

Finally, the group $T$ is maximal compact in $G$ and any compact subgroup can be
conjugated into $T$.  \end{Thm}

\begin{Rem}
This decomposition of $G=ANK$ is called an \emph{Iwasawa decomposition}. \end{Rem}

\begin{Rem}
By definition, a semisimple Lie algebra has no center.  Since a semisimple Lie group is
one for which its Lie algebra is semisimple, it follows that the center of a semisimple
Lie group is necessarily discrete, and that the $K$ in a Iwasawa decomposition is compact
if and only if the center is finite.  \end{Rem}

\begin{Cor} \label{IwasawaSplit}
Let $G$ be a semisimple Lie group with finite center, and $ANK$ be a Iwasawa
decomposition.  Then $G$ is bilipschitz diffeomorphic to $AN \times K$. \end{Cor}

\begin{pf}
The map $G \rightarrow AN \times K$ sends an element $p = ank \mapsto (an, k)$, and we
only need to check that the ratio between $d(a_{1}n_{1}k_{1}, a_{2}n_{2}k_{2})$ and
$d((a_{1}n_{1},k_{1}), (a_{2}n_{2},k_{2}))$ is bounded. This will follow from the
following calculation. Take $k_{i} \in K$, $g_{i} \in AN$,
\[ d(g_{1}k_{1}, g_{2}k_{2})=d(e, k_{1}^{-1}g_{1}^{-1}g_{2}k_{2})= d(e, k_{1}^{-1}k_{2}
\left( k_{2}^{-1} g_{1}^{-1}g_{2} k_{2} \right) ) \]

\noindent As $K$ is compact, $$ \max \{ \frac{d(e,g)}{d(e,kgk^{-1})},
\frac{d(e,kgk^{-1})}{d(e,g)}, g \in G, k \in K \} < \infty.$$
\end{pf}
Notice that when the center of $G$ is finite then it is easy to compute the
Assouad-Nagata dimension applying the previous corollary and the fact that $K$ is compact
in such case. The idea of the proof for  semisimple Lie groups (not necessarily with
finite center) will be to reduce the general case to the finite center case using
universal covers. So we have to focus in the fundamental group of a Lie group and its
universal cover.

\begin{Prop} \cite{Mil} (\textit{Proposition 1.6.4})
The fundamental group of any Lie group is a subgroup of its center. \end{Prop}

\begin{pf}
Let $G$ be a Lie group and $\tilde{G}$ its universal cover.  Then $\pi_{1}(G)$ is a
normal discrete subgroup of $\tilde{G}$ and $G=\tilde{G}/\pi_{1}(G)$.  Since $\pi_{1}(G)$
is normal, for a fixed $d \in \pi_{1}(G)$ we can define a map $\alpha: \tilde{G}
\rightarrow \pi_{1}(G)$ by $g \mapsto gdg^{-1}$.  Since this is a continuous map and
$\pi_{1}(G)$ is discrete, it follows that the image must be a point, namely $d$.  So this
shows that $\pi_{G}$ lies in the center of $\tilde{G}$. \end{pf}

\begin{Lem}
Let $G$ be a semisimple Lie group, and $ANK$ a Iwasawa decomposition.  Then $\tilde{G}$,
the universal cover of $G$, is also semisimple and $AN \tilde{K}$ is a Iwasawa
decomposition of $\tilde{G}$ where $\tilde{K}$ is the universal cover of $K$. \end{Lem}

\begin{pf}
That $\tilde{G}$ is semisimple follows from the fact that a covering map is a local
diffeomorphism and a Lie group is defined to be semisimple if its Lie algebra is
semisimple.  Since $A$ and $N$ are simply connected, it follows that
$\pi_{1}(G)=\pi_{1}(K)$, so the second claim follows.  \end{pf}

The next lemma will be used also in the next section. It is just a technical lemma about
abelian Lie groups. It will help to compute the Assouad-Nagata dimension of a group  of
the form $G/\Gamma$ with $G$ an abelian Lie group and $\Gamma$ a discrete subgroup.

\begin{Lem}\label{AbelianLemma}
Let $\Gamma$ be a discrete subgroup of $\ZZ^{k} \times \RR^{l}$.  Then $\left(  \ZZ^{k}
\times \RR^{l} \right)/ \Gamma$ is quasi-isometric to $\ZZ^{j}$ for some $j$.  \end{Lem}

\begin{pf}
As $\ZZ^{k} \times \RR^{l}$ is q.i. to $\ZZ^{k} \times \ZZ^{l}$, suffice to describe
$\left( \ZZ^{k} \times \ZZ^{l} \right)/ \Gamma$.  Since $\Gamma$ is a graph of a
homomorphism from $\ZZ^{\min \{ k,l\}} \rightarrow \ZZ^{\max \{k, l \}}$, by expressing
$\left( \ZZ^{k} \times \ZZ^{l} \right)$ in terms of a basis that contains a basis
representing $\Gamma$ as a graph, the quotient $\left( \ZZ^{k} \times \ZZ^{l} \right)/
\Gamma= \ZZ^{j}$ for some $j$.  \end{pf}

Now we provide the main theorem of this section.

\begin{Thm}\label{MainTheorem2}
Let $G$ ba a semisimple Lie group with a Riemannian metric $d$.  Then $asdim_{AN}(G, d) =
dim(G/T)$ where $T$ is a maximal compact subgroup of $G$  \end{Thm}

\begin{pf}
By corollary 3.6 of \cite{Car} we have $dim(G/T) = asdim(G) \le asdim_{AN}(G)$. So the
unique thing we have to prove is $asdim_{AN}(G) \le dim(G/T)$. \smallskip

By Theorem \ref{Iwasawa}, $G=ANK$, where $K= V \rtimes T$, and $V$ is connected, simply
connected abelian Lie group and $T$ is a maximal compact subgroup.  We cannot apply
Proposition \ref{ExactSequence} because $AN$ is not normal in $G$.  Instead, we observe
that $G/Z(G)$ is a semisimple Lie group with trivial center, and since $Z(G) < K$, it
follows that $\hat{G}=AN(K/Z(G))$ is a Iwasawa decomposition for $\hat{G}$, where
$K/Z(G)$ is compact.  By Corollary \ref{IwasawaSplit}, we have a bilipschitz
diffeomorphism $\hat{\phi}$ between $\hat{G}$ and $AN \times K/Z(G)$. Lifting
$\hat{\phi}$ up to $G$ we have a bilipschitz diffeomorphism $\phi$ between $G$ and $AN
\times K$ following commutative diagram.
\begin{equation*}
\begin{CD}
G  @>{\phi} >>  AN \times K \\
@VVV                                @VVV \\
\hat{G}                   @> {\hat{\phi}}>> AN \times K/Z(G)
\end{CD} \end{equation*}

So now $asdim_{AN}(G) \leq asdim_{AN}(AN) + asdim_{AN}(K)$.  But Theorem \ref{Iwasawa}
says that $K/V$ is compact and $V$ is isometrically embedded in $K$, it follows that
$asdim_{AN}(K)=asdim_{AN}(V)$, and so $asdim_{AN}(G) \leq asdim_{AN}(AN) + asdim_{AN}(V)
= dim(AN) + dim(V)$, since $AN$ and $V$ are both connected, simply connected solvable Lie
groups.  But $dim(AN) + dim(V)=dim(G/T)$.   \end{pf}

\section{Assouad-Nagata dimension of connected Lie groups}\label{MainSection3}
In this section we will prove the main result of this paper. We will reduce the general
case to the computations made for nilpotent, solvable and semisimple Lie groups in the
previous sections. In some sense the proof is similar to the one of solvable Lie group.
We will need to study an exact sequence generated by some subgroup  that is exponentially
distorted in our Lie group. Hence we will need the following concept from \cite{DeCor}:

\begin{Def}
Two Lie groups are said to be locally isomorphic if they have isomorphic Lie algebras.
\end{Def}

\begin{Rem} \label{WhatLocalIsomorphMeans}
If $G_{1}$, $G_{2}$ are two Lie groups with isomorphic Lie algebras, then there is a
simply connected Lie group $\tilde{G}$ and discrete subgroups $\Gamma_{1}, \Gamma_{2} <
Z(\tilde{G})$ such that $G_{1}=\tilde{G}/\Gamma_{1}$, $G_{2} = \tilde{G}/\Gamma_{2}$.
\end{Rem}

\begin{Def}If $G$ is a connected Lie group. Its exponential radical $R_{exp}(G)$ is the subgroup of $G$ so that
$G/R_{exp}(G)$ is the biggest quotient locally isomorphic to a direct product of a semisimple group and a group with
polynomial growth.
\end{Def}

\begin{Lem} \label{CorLem} Let $G$ be a connected Lie group and let $R_{exp}(G)$ be its exponential radical.
Then $R_{exp}(G)$ is strictly exponentially distorted in $G$ and it is contained in the nilpotent radical of $G$.
\end{Lem}
\begin{pf}
See Theorem 6.3. of \cite{DeCor}
\end{pf}

So now we need to study the quotient $G/R_{exp}(G)$. The key will be to understand what
means to be locally isomorphic from a large scale point of view.  For such purpose we
need the following:

\begin{Lem}\label{SimpleCenter} Let $G$ be simply connected semisimple Lie group.  Then $G$ is
quasi-isometric to $G/Z(G) \times Z(G)$.  \end{Lem}

\begin{pf}We already know that $G$ is quasi-isometric to $G/\tilde{K} \times \tilde{K}$ where
$\tilde{K}$ is the appropriate factor in the Iwasawa decomposition of $G$. Since $Z(G) <
\tilde{K}$ is a co-compact subgroup, the following sequence
$$ 1 \rightarrow G/ Z(G) \rightarrow G/\tilde{K} \rightarrow
\tilde{K}/Z(G) $$ shows that $G/Z(G)$ is quasi-isometric to $G/\tilde{K}$, to which
the desired claim follows.  \end{pf}

\begin{Thm}( \textit{Theorem 2.3 from \cite{Hochschild}}\label{ComopactAbelianSplit})
\newline
Leg $G$ be a topological group containing a vector group $V$ as a closed normal subgroup.
If $G/V$ is compact, then $V$ is a semidirect product of $G$.  \end{Thm}

\begin{Cor} \label{CompactSolvableSplit}
Let $G$ be a Lie group and $S  < G$ a normal solvable Lie group.  If $G/S$ is compact,
then $G$ is a semidirect product between a solvable group with the same topological
dimension as $S$ and a compact group. \end{Cor}

\begin{pf}
We induct on the length of the commutator series of $S$.  The base case is when $S$ is
abelian and this is given by Theorem \ref{ComopactAbelianSplit}.  Now suppose this is
true for solvable groups of length $j-1$.  Since $S$ is normal in $G$, it follows that
$Z(S)$, is also normal in $G$.  Therefore we can apply the inductive hypothesis to the
pair $G/Z(S)$ with subgroup $S/Z(S)$ and conclude that $G/Z(S)= \hat{S} \rtimes \hat{T}$,
where $\hat{S}$ is a solvable group and $\hat{T}$ a compact group.  Write
$\hat{S}=\tilde{S}/Z(S)$, and $\hat{T}=A/Z(S)$ for subgroups $\tilde{S}, A < G$.  Then
$G=\tilde{S} \rtimes A$.  We can apply Theorem \ref{ComopactAbelianSplit} to $A$ and
$Z(S)$ to conclude that $A= Z(S) \rtimes T$ for some compact subgroup $T$.  Therefore
$G=\tilde{S} \rtimes (Z(S) \rtimes T)$. \end{pf}

\begin{Rem}\label{SolvableCenterDistort}
Note that the center of a connected, simply connected solvable Lie group is at most
exponentially distorted.  To justify this claim, we need two observations. First, for
such a solvable Lie group $S$, the exponential map is a diffeomorphism and we can
coordinatize $S$ by its Lie algebra.  We now describe the second observation which is
taken from Section 2 of \cite{Breuillard}.  Given $\mathfrak{s}$ a solvable Lie algebra,
there exist a vector subspace $\mathfrak{v}$ and a nilpotent ideal $\mathfrak{n}$ (the
nilpotent radical) such that as a vector space, $\mathfrak{s} = \mathfrak{v} \oplus
\mathfrak{n}$, and for each $x \in \mathfrak{v}$, $\mathfrak{v}$ lies in the zero
generalized eigenspace of $ad(x)$. So if the Lie algebra of a subgroup lies in
$\mathfrak{v}$, then the distance of the subgroup is at most polynomially distorted. On
the other hand, since $\mathfrak{n}$ is a nilpotent ideal, whenever the Lie algebra of a
subgroup lies in $\mathfrak{n}$, the distance of the subgroup is at most
exponential-polynomially distorted. \smallskip

The claim now follows since we can always split the center of $S$ into a direct sum of
subgroups, whose Lie algebras lie in $\mathfrak{v}$ and $\mathfrak{n}$ respectively.
\end{Rem}

\begin{Thm}\label{MainTheorem3} Let $G$ be a connected Lie group.
Then $asdim_{AN}(G)= dim(G/C')$, where $C'$ is a maximal compact subgroup of $G$.
\end{Thm}

\begin{pf}
By Corollary 3.6 of \cite{Car} we have $dim(G/C) = asdim(G) \le asdim_{AN}(G)$, so
suffice to show that $asdim_{AN}(G) \le dim(G/T)$.\\

Let $R_{exp}(G)$ be the exponential radical of $G$.  Then we have the following short
exact sequence. \begin{equation}\label{ExpRadicalSSE} 1 \rightarrow R_{exp}(G)
\rightarrow G \rightarrow G/R_{exp}(G) \rightarrow 1,\end{equation} where quotient $U=
G/R_{exp}(G)$ is locally isomorphic to a direct product of a semisimple Lie group and a
solvable Lie group.  By Remark \ref{WhatLocalIsomorphMeans} there are simply connected
semisimple Lie group $L$ and a simply connected solvable Lie group $S$ such that $U$ is
isomorphic to $(L \times S)/\Gamma$ for some discrete subgroup $\Gamma < Z(L\times
S)=Z(L) \times Z(S)$. By Lemma \ref{SimpleCenter} we have
\begin{equation}\label{LocalIsoSSE} U \stackrel{q.i.}{\sim} L/Z(L) \times \left( Z(L)
\times S  \right) /\Gamma .\end{equation} In addition, we also have
\begin{equation}\label{CenterSSE} 1 \rightarrow \left(  Z(L) \times Z(S) \right) /\Gamma
\rightarrow \left( Z(L) \times S  \right) / \Gamma \rightarrow S/Z(S) \rightarrow 1
,\end{equation} \noindent where we note that $Z(S)$ is at most exponentially distorted in
$S$ by Remark \ref{SolvableCenterDistort}, and the same is true for $\left( Z(L) \times
Z(S) \right)/\Gamma$ in $\left( Z(L) \times S \right)/\Gamma$. \\

Putting (\ref{ExpRadicalSSE}), (\ref{LocalIsoSSE}) and (\ref{CenterSSE}) together we see
that \begin{eqnarray*}
asdim_{AN}(G) & \leq & asdim_{AN}(R_{exp}(G)) + asdim_{AN}(G/R_{exp}(G)) \\
& =& dim(R_{exp}(G)) + asdim_{AN}(U) \\
& \leq & dim(R_{exp}(G)) + asdim_{AN}(L/Z(L)) + asdim_{AN}(\left( Z(L) \times S \right)/\Gamma) \\
& = & dim(R_{exp}(G)) + dim(L/Z(L)) + asdim_{AN}(\left( Z(L) \times S \right)/\Gamma) \\
& \leq & dim(R_{exp}(G)) + dim(L/Z(L)) \\
&& + asdim_{AN}( \left( Z(L) \times Z(S) \right)/\Gamma) + asdim_{AN}( S/Z(S))  \\
& =&  dim(R_{exp}(G)) + dim(L/Z(L)) + dim( \left( Z(L) \times Z(S) \right)/\Gamma) + dim(
S/Z(S))  \end{eqnarray*}

Now let $L=ANK$ be a Iwasawa decomposition, and $K= V \rtimes T$ as given by Theorem
\ref{Iwasawa}.  Then $dim(AN)=dim(L/Z(L))$, and $$U = AN \left( (V \times S)/\Gamma
 \rtimes T\right) .$$  But $(V \times S) / \Gamma$ is diffeomorphic to a product $E
\times T_{1}$ where $T_{1}$ is a maximal subgroup in $(V \times S)/\Gamma$, and $E$ is a
manifold diffeomorphic to $\RR^{n}$, where $n=dim(E)= dim(S/Z(S)) + dim((Z(L) \times
Z(S)) / \Gamma)$.  In other words, $U= G/R_{exp}(G)$ is diffeomorphic to $AN \times E
\times C$ where $C$ is the compact subgroup $T_{1} \rtimes T$.  By Corollary
\ref{CompactSolvableSplit} we see that the subgroup $\tilde{T} < G$ such that
$\tilde{T}/R_{exp}(G) = C$ is a semidirect product between a solvable group $S_{1}$ of
the same dimension as $R_{exp}(G)$ and a maximal subgroup $C' < G$.  Therefore
\begin{eqnarray*}
dim(G/C') &=& dim(AN) + dim(E) + dim(R_{exp}(G)) \\
&= &dim(L/Z(L)) + dim(S/Z(S)) + dim((Z(L) \times Z(S)) / \Gamma) + dim(R_{exp}(G)) \end{eqnarray*} \end{pf}

We now have an analogous result to Corollary 3.6 of \cite{Car}.
\begin{Cor}
Let $\Gamma$ be a cocompact lattice in a connected Lie group $G$.  Then
$$ asdim_{AN}(\Gamma) = dim(G/K), $$ where $K$ is a maximal compact subgroup.  \end{Cor}

\begin{pf}
Since $\Gamma$ is quasi-isometric to $G$, it follows from Theorem \ref{MainSection2} that
\[ asdim_{AN}(\Gamma) = asdim_{AN}(G) = dim(G/K) \] \end{pf}

For the global Assouad-Nagata dimension the result is the following:

\begin{Cor}\label{EqualityForLie}
Let $G$ be a connected Lie group equipped with a Riemannian metric. Then $dim_{AN} (G)= dim(G)$.
\end{Cor}

\begin{pf}
On one hand we have $dim_{AN}(G) = \max\{asdim_{AN} (G), cdim(G)\}$, by corollary \ref{CapDimLie} we get $cdim(G)\le dim(G)$. Combining this inequality with the previous lemma we get $dim_{AN}(G) \le dim(G)$. \par
On the other hand $dim (G) \le dim_{AN}(G)$ (see the original paper of Assouad \cite{Assouad}).
\end{pf}

The asymptotic Assouad-Nagata dimension gives us an obstruction for a finitely generated
group to be quasi-isometrically embeddable in a connected Lie group in particular
quasi-isometric to a lattice. For example for certain classes of wreath products we have
the following:

\begin{Cor}
Let $H$ be a finite group ($H\ne 1$) and let $G$ be a finitely generated group such that
the growth is not bounded by a linear function. Then $H\wr G$ equipped with a word metric
can not be embedded quasi-isometrically in any connected Lie group.
\end{Cor}
\begin{pf}
In corollary 5.2 of \cite{Brod-Dydak-Lang} it was proved that $asdim_{AN}(H\wr G) =
\infty$. By Theorem \ref{MainTheorem3} the asymptotic Assouad-Nagata dimension of any
connected Lie group is finite. Therefore $H\wr G$ can not be embedded quasi-isometrically
in any connected Lie group.
\end{pf}

\begin{Rem}
Notice that the previous corollary implies that the solvable groups $\ZZ_2 \wr \ZZ^n$
with $n>1$ can not be a quasi-isometric to a cocompact lattice in a connected Lie group.
\end{Rem}

For general Lie groups we also have the following embedding result that is a consequence
of the results of \cite{Lang}. For more direct consequences we recommend to check such
paper.

\begin{Cor}Let $G$ be a connected Lie group equipped with a Riemannian
metric $d_{G}$.  Then there exists a $0 < p \le 1$ such that $(G, d_{G}^{p})$ can be
bilipschitz embedded into a product of $dim(G)+1$ many trees.
\end{Cor}
\begin{pf}
Combine \ref{EqualityForLie} with Theorem 1.3 of \cite{Lang}.
\end{pf}

\end{document}